\documentclass[11pt]{article}
\usepackage{amssymb,amsfonts,amsmath,latexsym,epsf,tikz,url}
\newtheorem{theorem}{Theorem}[section]

\newtheorem{conjecture}[theorem]{Conjecture}
\newtheorem{corollary}[theorem]{Corollary}

\newtheorem{definition}[theorem]{Definition}

\newcommand{\proof}{\noindent{\bf Proof.\ }}
\newcommand{\qed}{\hfill $\square$\medskip}

\begin{document}

\title{ More on the total dominator chromatic number of a graph}

\author{Nima Ghanbari and Saeid Alikhani$^{}$\footnote{Corresponding author} }

\date{\today}

\maketitle

\begin{center}

   Department of Mathematics, Yazd University, 89195-741, Yazd, Iran\\
{\tt n.ghanbari.math@gmail.com, alikhani@yazd.ac.ir}\\

\end{center}


\begin{abstract}

Let $G$ be a simple graph. A total dominator coloring of $G$, is a proper coloring of the vertices
of $G$ in which each vertex of the graph is adjacent to every vertex of some color class.
The total dominator chromatic (TDC) number $\chi_d^t(G)$ of $G$, is the minimum number of colors
among all total dominator coloring of $G$. The neighbourhood corona of two graphs $G_1$ and $G_2$ is denoted by  $G_1 \star G_2$  and is the graph obtained by
taking one copy of $G_1$ and  $|V(G_1)|$ copies of $G_2$, and joining the neighbours of the $i$th vertex of $G_1$ to every vertex in the $i$th copy of $G_2$. 
In this paper, we study the total dominator chromatic number of the neighbourhood of two graphs and  investigate the total dominator chromatic number of $r$-gluing of two  graphs.
Stability (bondage number) of total dominator chromatic number of $G$ is the minimum number of
vertices (edges) of $G$ whose removal changes the TDC-number of $G$. We study the stability and bondage number of certatin graphs. 

\end{abstract}


\medskip
\noindent{\bf AMS Subj.\ Class.:} 05C15, 05C69

\section{Introduction}

In this paper, we are concerned with simple finite graphs. Let $G=(V,E)$ be such a graph and $\lambda \in \mathbb{N}$. A mapping $f : V (G)\longrightarrow \{1, 2,...,\lambda\}$ is
called a $\lambda$-proper  coloring of $G$, if $f(u) \neq f(v)$ whenever the vertices $u$ and $v$ are adjacent
in $G$. A color class of this coloring, is a set consisting of all those vertices
assigned the same color. If $f$ is a proper coloring of $G$ with the coloring classes $V_1, V_2,..., V_{\lambda}$ such
that every vertex in $V_i$ has color $i$, then sometimes write simply $f = (V_1,V_2,...,V_{\lambda})$.  The chromatic number $\chi(G)$ of $G$ is
the minimum number of colors needed in a proper coloring of a graph.
The concept of a graph coloring and chromatic number is very well-studied in graph theory.

Vijayalekshmi and Kazemi \cite{Adel,Adel2,Vij1,Vij2} studied the total dominator coloring, abbreviated TD-coloring. Let $G$ be a graph with no
isolated vertex, the total dominator coloring  is a proper coloring of $G$ in which each vertex of the graph is adjacent
to every vertex of some (other) color class. The total dominator chromatic number, abbreviated TDC-number, $\chi_d^t(G)$ of $G$ is the minimum number of color classes in a TD-coloring of $G$.   Computation of the TDC-number is NP-complete (\cite{Adel}).
The TD-chromatic number of some graphs, such as paths, cycles, wheels and the complement of paths and cycles has computed in
\cite{Adel}. Also Henning in \cite{GCOM} established the  lower and the upper bounds on the TDC-number
of a graph in terms of its total domination number $\gamma_t(G)$. He has shown that,  every
graph $G$ with no isolated vertex satisfies $\gamma_t(G) \leq \chi_d^t (G)\leq \gamma_t(G) + \chi(G)$.
The properties of TD-colorings in trees has studied in \cite{GCOM,Adel}. Trees $T$ with $\gamma_t(T) =\chi_d^t(T)$ has characterized
in \cite{GCOM}. We have examined the effects on $\chi_d^t(G)$ when $G$ is modified by operations on vertex and edge of $G$ in \cite{nima2}.  

\begin{figure}[ht]
	\begin{center}
		\includegraphics[width=0.37\textwidth]{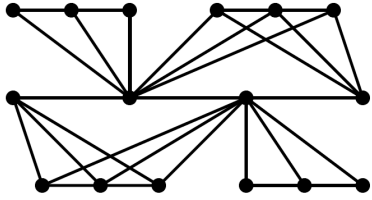}
		\caption{\label{fig1} The neighbourhood corona of $P_4 \star P_3$.}
	\end{center}
\end{figure}

The corona of two graphs $G$ and $H$ which denoted by $G\circ H$ is defined in \cite{Harary} and there have been some results on the
corona of two graphs \cite{Frucht}. In \cite{nima1} we have studied the total dominator chromatic number  of corona of two graphs. In this paper we consider another variation of corona of two graphs and study its total dominator chromatic number. 
Given simple graphs  $G_1$ and $G_2$, the neighbourhood corona of $G_1$ and $G_2$, denoted by $G_1 \star G_2$ and  is the graph obtained by
taking one copy of $G_1$ and $|V(G_1)|$ copies of $G_2$ and  joining the neighbours of the $i$th vertex of $G_1$ to every vertex in the $i$th copy of $G_2$ (\cite{Gopalapillai}). 
Figure \ref{fig1} shows $P_4\star P_3$, where $P_n$ is  the path of order $n$. 
Liu and Zhu in \cite{linear} determined the adjacency spectrum of $G_1\star G_2$ for arbitrary $G_1$ and $G_2$ and the Laplacian spectrum and signless Laplacian spectrum of $G_1\star G_2$ for regular $G_1$ and arbitrary $G_2$, in terms of the corresponding spectrum of $G_1$ and $G_2$. Also Gopalapillai in \cite{Gopalapillai} has studied 
the eigenvalues and spectrum of $G_1\star G_2$, when $G_2$ is regular. 

\medskip 

A domination-critical (domination-super critical, respectively) vertex in a graph
$G$, is a vertex whose removal decreases (increases, respectively) the domination
number. Bauer et al. \cite{1}  introduced the concept of domination stability in graphs.
The domination stability, or just $\gamma$-stability, of a graph $G$ is the minimum number
of vertices whose removal changes the domination number. Motivated by domination stability, we introduce the total dominator chromatic (TDC)-stability (TDC-bondage number) of a graph. 

\begin{definition}
The total dominator chromatic (TDC)-stability (TDC-bondage number) of graph  $G$,  is the minimum number of
vertices (edges) of $G$, whose removal changes the TDC-number of $G$.  
\end{definition}

\medskip 

In the next section, we study the total dominator chromatic number of neighbourhood corona of two graphs. We investigate the total dominator chromatic number of $r$-gluing of two graphs in Section 3. We study the TDC-stability and TDC-bondage number of certain graphs in Section 4.


\section{TDC-number of neighbourhood corona of two graphs}

The TDC-number of binary graph operations, aside from Cartesian  and   corona product (\cite{nima1,Adel2}), have not been widely studied. 
In this section, we shall study the total dominator chromatic (TDC) number of neighbourhood corona of two graphs. We begin with the following theorem which gives an upper bound for $\chi_d^t(G_1\star G_2)$. 
\begin{theorem}\label{1}
For any two connected graphs $G_1$ and $G_2$, 
	$$\chi_d^t(G_1\star G_2)\leq\vert V(G_1)\vert+\vert V(G_2)\vert.$$
\end{theorem}

\proof 
First we give the colors $a_1,a_2,\ldots,a_{|V(G_1)|}$ to the vertices of $G_1$. Now for every copy of $G_2$ we give the colors $b_1,b_2,\ldots,b_{|V(G_2)|}$ to its vertices. Since $G_1$ is connected, each vertex of $G_1$ uses one of its adjacent colors for TD-coloring. Also each vertex of each $G_2$ is adjacent  to at least one vertex of $G_1$ and can use this color for TD-coloring. Therefore this is a TD-coloring for the $G_1\star G_2$ and we have the result. \qed

\medskip 
The following theorem gives another upper bound for the total dominator chromatic number of the neighbourhood of two graphs.

\begin{theorem}\label{2}
	For any two connected graphs $G_1$ and $G_2$, 
	$$\chi_d^t(G_1\star G_2)\leq\chi_d^t(G_1)+\vert V(G_2)\vert.$$
\end{theorem}

\proof
First we give the colors $a_1,a_2,\ldots,a_{\chi_d^t(G_1)}$ to the vertices of $G_1$ as we colored  $G_1$ for a TD-coloring.   Now for every copy of $G_2$ we give the colors $b_1,b_2,\ldots,b_{|V(G_2)|}$ to its vertices. Now we consider one copy of $G_2$ and we call  the corresponding vertex for this $G_2$, say $w$. Each vertex of $G_2$ is now adjacent to all of the adjacent vertices of $w$. So each color of $w$ for the TD-coloring, can use for each vertices of $G_2$, too. So this is a TD-coloring of $G_1\star G_2$ and we have 
$$\chi_d^t(G_1\star G_2)\leq\chi_d^t(G_1)+\vert V(G_2)\vert.$$
\qed


\begin{theorem}\label{3}
	For any two connected graphs $G_1$ and $G_2$, 
	$$\chi_d^t(G_1\star G_2)\leq\chi_d^t(G_1)+\chi_d^t(G_2).$$
\end{theorem}

\proof
We give the colors $a_1,a_2,\ldots,a_{\chi_d^t(G_1)}$ to the vertices of $G_1$ as we colored $G_1$ for a TD-coloring. So each vertex of $G_1$ uses the old color class. Now for every copy of $G_2$ we give the colors $b_1,b_2,\ldots,b_{\chi_d^t(G_2)}$ to its vertices as we color $G_2$ for a TD-coloring. Now by the same argument as the proof of the Theorem \ref{2} we have the result.
\qed 

The following theorem present an upper bound for $\chi_d^t(G_1\star G_2)$ based on total dominator chromatic number of $G_1$ and chromatic number of $G_2$. 

\begin{theorem}\label{4}
	For any two connected graphs $G_1$ and $G_2$, 
	$$\chi_d^t(G_1\star G_2)=\chi_d^t(G_1)+\chi(G_2).$$
\end{theorem}

\proof
We give the colors $a_1,a_2,\ldots,a_{\chi_d^t(G_1)}$ to the vertices of $G_1$ as we colored $G_1$ for a TD-coloring. So each vertex of $G_1$ uses the old color class. Now for every copy of $G_2$ we give the colors $b_1,b_2,\ldots,b_{\chi(G_2)}$ to its vertices as we color $G_2$ for a coloring. Now by the same argument as the Proof of the Theorem \ref{2} we have: 
$$\chi_d^t(G_1\star G_2)\leq\chi_d^t(G_1)+\chi(G_2).$$
Now we consider one color class such as $j\in\{a_1,a_2,\ldots,a_{\chi_d^t(G_1)},b_1,b_2,\ldots,b_{\chi(G_2)}\}$. we cannot omit the color class $j$, because the coloring is proper and all vertices of $G_2$ are adjacent to at least one vertex of $G_1$ and since $G_1$ is connected, every vertex of $G_1$ is adjacent to at least all of the vertices of a $G_2$. So we have the result.
\qed

\noindent{\bf Remark 1.} The bounds in the Theorems \ref{1},\ref{2}, and \ref{3} are sharp. It suffices to consider the complete graph $K_4$ as $G_1$ and $K_3$ as $G_2$.

\medskip

The friendship graph $F_n$  $(n \geq  2)$ can be constructed by joining $n$ copies of the cycle graph $C_3$ with a common vertex. It is easy to see that $\chi_d^t(F_n)=3$ (\cite{nima1}). We end this section with the following corollary which follows from Theorem \ref{4}:

\begin{corollary}
	\begin{enumerate}
		\item [(i)] For every natural number $n$, $\chi_d^t(F_n\star K_n)=n+3$.
		
		\item[(ii)] For every natural number $n$, $\chi_d^t(F_n\star C_{2n})=5.$
		\end{enumerate}
		\end{corollary} 
		
\section{TDC-number of $r$-gluing of two graphs}
Let $G_1$ and $G_2$ be two graphs and $r\in \mathbb{N}\cup \{0\}$ with $r\leq min \{\omega(G_1),\omega(G_2)\}$, where $\omega(G)$ shows the clique number of $G$.  Choose a $K_r$ from each $G_i$, $i=1,2$, and form a new graph $G$ from the union of $G_1$ and $G_2$ by identifying the two chosen $K_r$'s in an arbitrary manners. The graph $G$ is called $r$-gluing of $G_1$ and $G_2$ and denoted by $G_1\cup_{K_r}G_2$. If $r=0$ then $G_1\cup_{K_0}G_2$ is just disjoint union. The $G_1\cup_{K_i}G_2$ for $i=1,2$, is called vertex and edge gluing, respectively. 
In this section, we study the total dominator chromatic number of $r$-gluing of two graphs.

\begin{theorem}\label{glue}
	For any two connected graphs $G_1$ and $G_2$, 
	$$max \{\chi_d^t(G_1),\chi_d^t(G_2)\}\leq\chi_d^t(G_1\cup_{K_r} G_2)\leq\chi_d^t(G_1)+\chi_d^t(G_2)-r.$$
\end{theorem}

\proof
Since we need at least $\chi_d^t(G_1)$ colors to color $G_1$ and $\chi_d^t(G_2)$ colors to color $G_2$, so we need at least max $\{\chi_d^t(G_1),\chi_d^t(G_2)\}$ colors to color $G_1\cup_{K_r} G_2$. So we have $max \{\chi_d^t(G_1),\chi_d^t(G_2)\}\leq\chi_d^t(G_1\cup_{K_r} G_2)$.

On the other hand, first we can give $a_1,a_2,\ldots,a_r$ to the vertices of $K_r$. Then we give $a_{r+1},\ldots,a_{\chi_d^t(G_1)}$ to the other vertices of $G_1$ to have a TD-coloring for $G_1$. Also we give $b_1,b_2,\ldots,b_{\chi_d^t(G_1)-r}$ to the other vertices of $G_2$ to have a TD-coloring for $G_2$. So every vertex of $G_1\cup_{K_r} G_2$ uses the color class which used before and this is a TD-coloring for $G_1\cup_{K_r} G_2$. So $\chi_d^t(G_1\cup_{K_r} G_2)\leq\chi_d^t(G_1)+\chi_d^t(G_2)-r.$
\qed

\medskip
\noindent{\bf Remark 2.} The bounds in the Theorem \ref{glue} are sharp. For the lower bound it suffices to consider the complete graph $K_4$ as $G_1$ and $K_5$ as $G_2$ and $r=4$. For the upper bound it suffices to consider the cycle graph $C_4$ as $G_1$ and $K_3$ as $G_2$ and $r=1$.

\section{TDC-stability and TDC-bondage number of  certain graphs}

In this section, we study the stability and bondage number of total dominator chromatic number of certain graphs. First we consider stability of certain graphs.  

\subsection{TDC-stability of certain graphs}

Stability of total dominator chromatic number of a graph $G$, $St_d^t(G)$,  is the minimum number of vertices of $G$و whose removal changes the TDC-number of $G$. To obtain the 
stability of specific graphs, we need the following results:  
\begin{theorem}{\rm(\cite{Adel})}\label{CnPn}
	\item[(i)]
	Let $P_n$ be a path of order $n\geq 2$. Then
	\[
	\chi_d^t(P_n)=\left\{
	\begin{array}{lr}
	{\displaystyle
		2\lceil\frac{n}{3}\rceil-1},&
	\quad\mbox{if $n\equiv 1$ $(mod\,3)$,}\\[15pt]
	{\displaystyle
		2\lceil\frac{n}{3}\rceil},&
	\quad\mbox{otherwise.}
	\end{array}
	\right.
	\]
	\item[(ii)] Let $C_n$ be a cycle of order $n\geq 4$. Then
	\[
	\chi_d^t(C_n)=\left\{
	\begin{array}{lr}
	{\displaystyle
		2,} &
	\quad\mbox{if $n=4$}\\[15pt]
	{\displaystyle
		4\lfloor\frac{n}{6}\rfloor+r,}&
	\quad\mbox{if $n\neq 4$, $n\equiv r$ $(mod\,6)$, $r=0,1,2,4$,}\\[15pt]
	{\displaystyle
		4\lfloor\frac{n}{6}\rfloor+r-1,}&
	\quad\mbox{if $n\equiv r$ $(mod\,6)$, $r=3,5$.}
	\end{array}
	\right.
	\]
\end{theorem}

\begin{theorem}\label{Pn}
 For any  $n\geq 4$,  $St_d^t(P_n)=1$. 
\end{theorem}

\begin{figure}
		\begin{center}
			\includegraphics[width=5in]{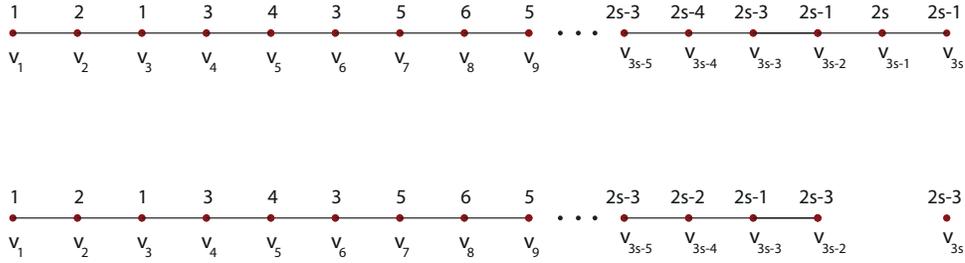}
			\caption{Coloring of $P_{3s}$ and removing the vertex $v_{3s-1}$. }
			\label{P3s}
		\end{center}
	\end{figure}

\proof
We prove the theorem for  three following cases:
\begin{itemize}
\item[(i)] If $n\equiv 0$  $(mod\,3)$, then we have shown a TDC coloring of $P_{3s}$ in Figure \ref{P3s}. By removing the vertex  $v_{3s-1}$,  we have a TDC coloring  by $3s-1$ colors. So we have $St_d^t(P_{3s})=1$.

\item[(ii)] If $n\equiv 1$  $(mod\,3)$, then in this case we remove an end  vertex $v$. By Theorem \ref{CnPn}, we know that $\chi_d^t(P_{3s+1})=2s+1$ and $\chi_d^t(P_{3s})=2s$. So we have $St_d^t(P_{3s+1})=1$.

\item[(iii)] If $n\equiv 2$  $(mod\,3)$, the proof is similar to the proof of Part (ii).

\end{itemize}
\qed


\begin{theorem}\label{cycle}
 For every $n\geq 3$, 
$	{St}_d^t(C_n)=\left\{
	\begin{array}{lr}
	{\displaystyle
		1,} &
	\quad\mbox{if  $n\equiv r$ $(mod\,6)$, $r=1,2,4,5$,}\\[15pt]
	{\displaystyle
		1,}&
	\quad\mbox{if $n=3$}\\[15pt]
	{\displaystyle
		2,} &
	\quad\mbox{if $n\neq3$  $n\equiv r$ $(mod\,6)$, $r=0,3$.}
	\end{array}
\right.	$
\end{theorem}

\proof There are six cases, $n\equiv r$  $(mod\,6)$ and $r=0,1,2,3,4,5$. We only prove  one case and the proof of another  cases are similar. Suppose that $n=6k+3$. By Theorem \ref{CnPn}, we have $\chi_d^t(C_{6k+3})=4k+2$. By removing each vertex of this graph, we have a path graph of order $6k+2$ and By Theorem \ref{CnPn}, we have $\chi_d^t(P_{6k+2})=4k+2$. Now by Theorem \ref{Pn}, we have $St_d^t(P_{6k+2})=1$. Therefore we have $St_d^t(C_{6k+3})=2$.
\qed

\medskip
The friendship graph $F_n$ $(n\geqslant 2)$ can be constructed by joining $n$ copies of the cycle graph $C_3$ with a common vertex. 
The $n$-book graph $(n\geqslant 2)$ is defined as the Cartesian product $K_{1,n}\square P_2$. We call every $C_4$ in the book graph $B_n$, a page of $B_n$. All pages in $B_n$ have a common side $v_1v_2$. 
The following easy result gives the TDC-stability of $F_n$ and $B_n$.

\begin{theorem}
	\begin{enumerate}
		\item [(i)]
 For every  $n\geq 2$,  $St_d^t(F_n)=1$. 
 \item[(ii)] 
  For every $n\geq 3$, $St_d^t(B_n)=1$. 
  \end{enumerate} 
\end{theorem}

\begin{figure}
		\begin{center}
			\includegraphics[width=3.5in]{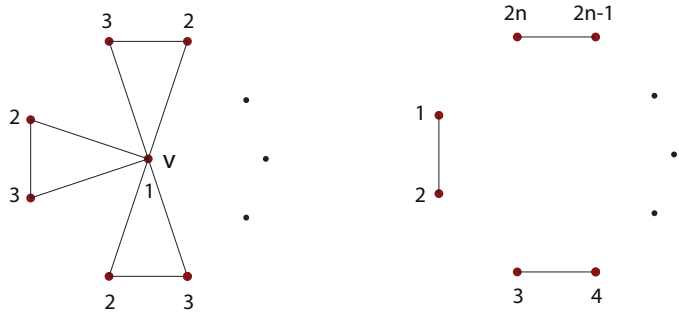}
			\caption{Coloring of $F_n$ and removing the vertex $v$. }
			\label{friend}
		\end{center}
	\end{figure}

\begin{figure}
		\begin{center}
			\includegraphics[width=3.5in]{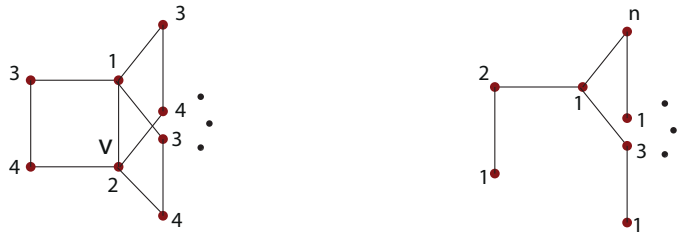}
			\caption{Coloring of $B_n$ and removing vertex $v$. }
			\label{book}
		\end{center}
	\end{figure}

\proof
	\begin{enumerate}
		\item [(i)]
As we see in the Figure \ref{friend}, by removing the  vertex  $v$ we have a TDC coloring for that graph using  $2n$ colors. So we have $St_d^t(F_n)=1$.
\item[(ii)] 
By removing the  vertex  $v$ in Figure \ref{book}, we have a TDC coloring for that graph using  $n$ colors. So we have $St_d^t(B_n)=1$.
\qed
 \end{enumerate}

\begin{theorem}
 For every $n\in \mathbb{N}$, there exists a graph $G$ such that  $ St_d^t(G)=n$.
\end{theorem}

\begin{figure}
		\begin{center}
			\includegraphics[width=1.6in]{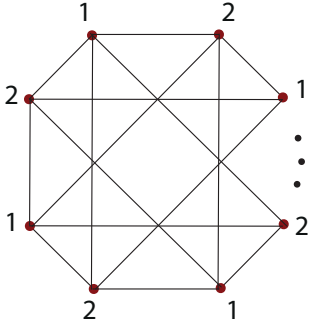}
			\caption{$n$-regular graph of order $2n$. }
			\label{nregular}
		\end{center}
	\end{figure}

\proof
Consider the graph $G$ of order $2n$ in Figure \ref{nregular}. As observe that,  each vertex with color $1$ is adjacent to every vertex  with color $2$ and each vertex with color $2$ is adjacent to every vertex  with color $1$, and so $\chi_d^t(G)=2$. By removing just one vertex of $G$, the coloring does not change. Suppose that $A$ is the set of vertices whose have  color $1$. The TDC-number of the induced graph $G-A$ is $1$. The set $A$  has  the minimum number of vertices which changes the TDC number  of these kind of graphs (since $K_2$ is always a subgraph of these graphs and we do not need to change the color of the graph by removing each vertex).  Therefore $ St_d^t(G)=n$.
\qed

	In 1956, Nordhaus and Gaddum obtained the lower and upper bounds for the sum
	of the chromatic numbers of a graph and its complement (actually, the upper bound
	was first proved by Zykov \cite{Zykov} in 1949). Since then, Nordhaus-Gaddum type bounds
	were obtained for many graph invariants. An exhaustive survey is given in \cite{Aouchiche}. So investigation of  Nordhaus-Gaddum type inequalities for the TDC-stability number is 
	an interesting problem and here we just present a lower bound and finding sharp upper bounds remain as an open problem.   
	Since we have to remove at least one vertex of both graph $G$ and $G^c$  to change  the TDC-number of $G$ and $G^c$, so we have the following easy result.

\begin{theorem}\label{GGc}
	For every graph $G$,  $St_d^t(G)+St_d^t(G^c)\geq 2$.
\end{theorem}

\noindent{\bf Remark 4.} There are graphs $G$ whose satisfy the equality in Theorem \ref{GGc}. Two graphs $P_3$ and $C_5$ are two examples.  
 As another example we can consider the graph $K_n-e$, where $e\in E(K_n)$.

	\medskip 
	
	We  end this subsection by proposing the following conjecture:

	\begin{conjecture} 
	Let $G$ be a graph with  a vertex of  degree one or two. Then $St_d^t(G)=1$ or $St_d^t(G)=2$. 
\end{conjecture} 

\subsection{TDC-bondage number number of certain graphs}

Bondage number of the total dominator coloring of a graph $G$, $B_d^t(G)$, is the minimum number of edges of $G$, whose removal changes the TDC-number of $G$.
In this subsection, we study the TDC-bondage number of specific graphs. 

\begin{theorem}\label{BPn}
 For the path graphs $P_n$ with $n\geq 3$,  $B_d^t(P_n)=1$.
\end{theorem}

\begin{figure}
		\begin{center}
			\includegraphics[width=5in]{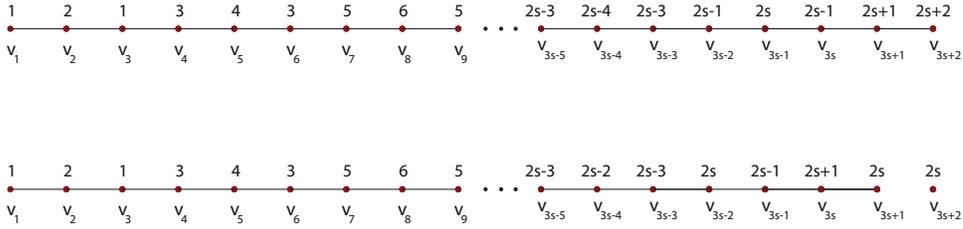}
			\caption{Coloring of $P_{3s+2}$ and removing  an edge. }
			\label{P3s2}
		\end{center}
	\end{figure}

\proof
We have three cases:
\begin{itemize}
\item[(i)] If $n\equiv 2$  $(mod\,3)$. In this case, as we see in the Figure \ref{P3s2}, by removing the edge between two vertices $v_{3s+1}$ and $v_{3s+2}$, we have a TDC coloring for new  graph by $3s+1$ colors. So  $B_d^t(P_{3s+2})=1$.

\item[(ii)] If $n\equiv 1$  $(mod\,3)$. In this case,  by removing the edge $v_2v_3$ we have  two paths $P_2$ and $P_{3s-1}$ and by Theorem \ref{CnPn},  $\chi_d^t(P_{3s+1})=2s+1$, $\chi_d^t(P_{3s-1})=2s$ and $\chi_d^t(P_2)=2$. So we have $B_d^t(P_{3s+1})=1$.

\item[(iii)] If $n\equiv 0$  $(mod\,3)$ and so $n=3s$ for some $s\in \mathbb{N}$. By removing a suitable edge  which makes two paths $P_{3t+2}$ and $P_{3k+1}$ ($t+k+1=s$), by Theorem \ref{CnPn}, we have $\chi_d^t(P_{3t+2})=2t+2$, and $\chi_d^t(P_{3k+1})=2k+1$. So $\chi_d^t(P_{3t+2})+\chi_d^t(P_{3k+1})=2t+2k+3$. On the other hand $\chi_d^t(P_{3s})=2s=2s+2k+2$. Therefore   $B_d^t(P_{3s+1})=1$.
\qed

\end{itemize}

\begin{theorem}
 For the cycle graphs $C_n$ with $n\geq 5$ we have:
 \[
	B_d^t(C_n)=\left\{
	\begin{array}{lr}	
   {\displaystyle
		1,} &
	\quad\mbox{if  $n\equiv 4$ $(mod\,6)$,}\\[15pt]
	{\displaystyle
		2,} &
	\quad\mbox{Otherwise.}
		\end{array}
	\right.
	\]
 
\end{theorem}

\proof
There exist  six cases, $n\equiv r$  $(mod\,6)$ and $r=0,1,2,3,4,5$. We only prove  one case and the proof of another cases are similar. Suppose that $n=6k+3$. By Theorem \ref{CnPn}, we have $\chi_d^t(C_{6k+3})=4k+2$. By removing each edge of this graph, we have a path graph of order $6k+3$ and by Theorem \ref{CnPn}, we have $\chi_d^t(P_{6k+3})=4k+2$. Now by Theorem \ref{BPn}, we have $B_d^t(P_{n})=1$. Therefore we need to remove two edges and so $B_d^t(C_{6k+3})=2$.
\qed

\begin{theorem}
 The bondage number of  the friendship graph $F_n$ ($n\geq 2$) is  $B_d^t(F_n)=1$.
\end{theorem}

\begin{figure}
		\begin{center}
			\includegraphics[width=3.5in]{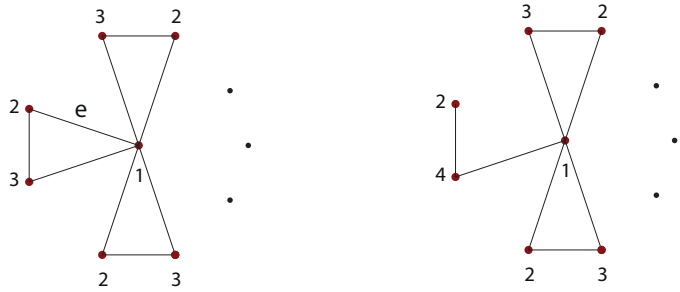}
			\caption{Coloring of $F_n$ and removing edge $e$. }
			\label{friend2}
		\end{center}
	\end{figure}

\proof
As we see in the Figure \ref{friend2}, by removing the  edge $e$ we have a TDC coloring for that graph by $4$ colors. So we have $B_d^t(F_n)=1$.
\qed

We end the paper by the following theorem and its remark:  
\begin{theorem}
 For every graph $G$,  $B_d^t(G)+B_d^t(G^c)\geq 2$.
\end{theorem}

\proof
Since for  changing  the TDC-number of $G$ and $G^c$,  we have to remove at least one edge of these two graphs, so  we have the result.
\qed

\begin{figure}
		\begin{center}
			\includegraphics[width=2.8in]{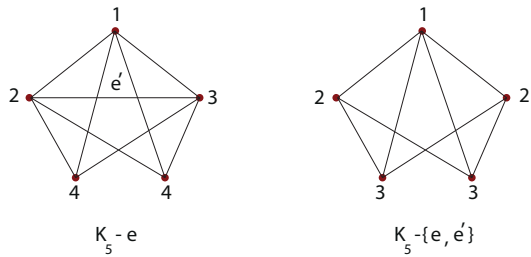}
			\caption{Coloring of $K_5-e$ and removing edge $e^\prime$. }
			\label{K5}
		\end{center}
	\end{figure}

\noindent{\bf Remark 5.} The lower bound in the Theorem \ref{GGc} is sharp. It is sufficient to consider Path graph $P_3$ as $G$. Also we can consider the graph $K_5-e$ as shown in Figure \ref{K5}.

\end{document}